\documentstyle [12pt,twoside, epsf]{article}

\def\picill#1by#2(#3)
{\vbox to #2
{\hrule width #1 height 0pt depth 0pt
\vfill\epsffile{#3}}}

\begin{document}

\title{\bf A Polynomial Invariant for Flat Virtual Links}

\author{Louis H. Kauffman \\
Department of Mathematics, Statistics and
Computer Science\\
University of Illinois at Chicago\\
851 South Morgan St., Chicago IL 60607-7045, USA\\
kauffman@uic.edu\\
and\\
R. Bruce Richter \\
Department of Combinatorics and Optimization\\
University of Waterloo\\
Waterloo, ON N2L 3G1, Canada\\
brichter@uwaterloo.ca
}

\date{December 2005}

\maketitle

\section{Introduction}

In searching for methods to distinguish knots and links, recent work has focussed on the more general virtual knots and links,
introduced by Kauffman in \cite{kauf1}.  If we ignore the over/under information at each crossing, we obtain a flat  virtual knot or link
(still called a knot or a link), sometimes
called the {\em universe\/} of the original virtual knot or link.  It is a triviality that if the original virtual link is the unlink, then the flat
virtual link is also the unlink.  Thus, being able to determine that a flat virtual link is not the unlink implies that every virtual link with that
universe is non-trivial.
\bigbreak

The focus of this work is on creating a polynomial invariant for flat virtual links.   
The polynomial is extremely simple to compute and, in the case of one component, specializes to  Turaev's virtual string polynomial \cite{turaev}. 
We shall see that Turaev's polynomial has the property that it is non-zero precisely when there is no filamentation of the knot, as described in 
Hrencecin and Kauffman \cite{hren}.  
\bigbreak

Schellhorn has provided a version of filamentations for flat virtual links \cite{schell}.  Our motivation was to try to turn
Schellhorn's filamentations into a polynomial invariant.  We have only partly succeeded:  our polynomial has the property that if there is a
filamentation, then the polynomial is 0.  Unfortunately the converse fails, although if the polynomial is 0, then it turns out to be easy to determine
if there is a filamentation.
\bigbreak

\noindent {\bf Acknowledgement.} LHK acknowledges support from NSF Grant DMS-0245588.
RBR acknowledges the financial support of NSERC.
\bigbreak

\section{The Polynomial}

Let $L$ be a flat virtual link.  There are three ways to represent $L$, illustrated in Figure 1.
(For a fuller discussion of the relationships between these representations, see, for example, \cite{hren} or \cite{schell}.).   The first
representation is the most geometric.  This has closed curves in the plane, with some crossings distinguished as ``virtual" (i.e, not ``really"
there).  These are the crossings inside a small circle in the diagram.   The second is the  {\em chord diagram\/} of $L$, which has a circle for each
component $A$ of $L$ and, for each crossing of $L$, there is a chord joining the two points on the circles representing that crossing.  The chord
has a direction, with the tail being the ``negative" crossing (the negative strand of the crossing has the other strand crossing from left to right) and
the head being the ``positive" crossing (the positive strand of the crossing has the other strand crossing from right to left).  The third
representation, which is the one we shall mainly use, is the {\em Gauss code\/}, which replaces the chords of the chord diagram with labels on the ends
of the chords.  A crossing $x$ will have the negative end of the chord labeled $x^-$ and the positive end labeled $x^+$.   The Gauss code is then simply
read off as the cyclic sequences of signed labels on each circle.
\smallbreak

There are no virtual crossings in the chord diagram or Gauss code representations of a virtual link. The virtual crossings in the planar diagram
are artifacts of this form of representation. Consequently, there is no need to consider the virtual crossings to prove invariance of definitions
that are made at the Gauss code or chord diagram level.
\smallbreak

If there is a pair $A,B$ of link components such that the number of $+$ ends of $AB$-chords on $A$ is not equal
to the number of $-$ ends of $AB$-chords on $A$, then the {\it flat linking number} of $A$ and $B$, which is just (half) the difference of the
number of $+$ ends on $A$ and the number of $-$ ends on $A$, is different from 0. So $A$ and $B$ are
linked, and therefore the entire virtual link is non-trivial.  Thus, we may assume that each pair of codewords $A$ and $B$ have an equal number of
positive and negative $AB$-letters on
$A$.  This implies that the sum of all the signs of letters in $A$ is 0.
\bigbreak

The set $\mathcal K(A)$ is the set of all self-crossings $x$ of $A$, i.e., both $x^+$ and $x^-$ are in the codeword for $A$.
For each unordered pair $\{A,B\}$ of components, we arbitrarily partition the $AB$-crossings into pairs $\{x,y\}$,
so that $x^+$ and $y^-$ are on the
same one of $A$ and $B$.  We denote by $\mathcal P_{A,B}$ the set of pairs in the partition.
\smallbreak

$$ \picill4inby2in(Figure1)  $$
\begin{center}
{\sc Figure 1 - A Flat Link}
\end{center}

Our polynomial $p_L$ will have a variable $t_A$ for each codeword $A$ of the link and a variable $t_{A,B}$
for each unordered pair of distinct codewords $A$ and $B$ of the link.
For two letters $x^\epsilon$ and $y^\delta$ of the same codeword $A$ of the link, define
$\eta_A(x^\epsilon,y^\delta)$ to be the difference of the number of $+$ ends  and the number of $-$ ends as we cycle through $A$ from $x^\epsilon$ to
$y^\delta$, not including these two letters.   The reader will have no difficulty realizing that this is the {\em intersection number\/}  of the
filament $x^\epsilon y^\delta$ (Definition 4.7 in \cite{hren} or Definition 3.5 in \cite{schell}).  See Figure 2.
\smallbreak

$$ \picill4inby2in(Figure2)  $$
\begin{center}
{\sc Figure 2 - A Chord Diagram}
\end{center}

The polynomial $p_L$ is the sum of two polynomials.  The first polynomial is the sum, over all codewords $A$ of the link, of the polynomials
$$p_A(t_A)=\sum_{x\in\mathcal K(A)}\eta_A(x^+,x^-)t_A^{|\eta_A(x^+,x^-)|}\,.$$
The reader should recognize $p_A(t_A)$ as Turaev's  virtual string polynomial, $u(t_{A})$ \cite{turaev}.
\bigbreak

The second polynomial is the sum over all pairs $\{A,B\}$ of distinct components $A$ and $B$ of the link of the linear term
$$p_{A,B}(t_{A,B})=t_{A,B}\sum_{\{x,y\}\in \mathcal P(A,B)}\left(\eta_A(x^+,y^-)+\eta_B(y^+,x^-)\right)\,,$$ 
where the labelling of $x$ and $y$ is always to be chosen so that $x^+$ and $y^-$ are on $A$, while $x^-$ and $y^+$ are on $B$.
\bigbreak

\noindent Our main theorem is:
\smallbreak

\noindent {\bf Theorem A.} {\it The polynomial $p_L$ is an invariant of the flat virtual link $L.$}
\bigbreak

The proof is broken up into two steps:  invariance of $p_{A,B}$ relative to
the partition $\mathcal P_{A,B}$ and invariance relative to the flat virtual Reidemeister moves.  The first of these is taken up in the next section,
while the Reidemeister moves are taken up in Section \ref{ReidInv}.  As mentioned earlier, we may ignore those flat virtual Reidemeister moves that
involve virtual crossings.
\bigbreak

\section{Invariance relative to partitions}
In this section, we show that $p_{A,B}$ is independent of the particular partition
${\mathcal P}(A,B)$.  First observe that any partition is obtained from a particular one by a finite sequence of {\em elementary switches\/}, i.e.,
replacing two pairs $\{e,f\}$ and $\{e',f'\}$ by the two pairs $\{e,f'\}$ and $\{e',f\}$, where $e$ and $e'$ both have their positive ends on $A$.  So
it is sufficient to show that $p_{A,B}$ is unaffected by an elementary switch.
\bigbreak

In the codeword $A$, there are only two possible cyclic orderings of the signs of
the four points $e^+$, $e'{}^+$, $f^-$, and $f'{}^-$, namely $++--$ and $+-+-$.  For each of these, there are two ways to pair each $+$ with a $-$.  If
we have $a^+,b^+,c^-,d^-$, then the pairings $\{a,c\}$ and $\{b,d\}$ count each $+/-$ between (never including endpoints) $a^+$ and $b^+$ once, between
$b^+$ and $c^-$ twice, between $c^-$ and $d^-$ once, and between $d^-$ and $a^+$ not at all.  In addition, the $+$ in $b^+$ and the $-$ in $c^-$ are
both counted once.  On the other hand, in the pairings $\{a,d\}$ and $\{b,c\}$, we have exactly the same counts.  Thus,
$\eta_A(a^+,c^-)+\eta_A(b^+,d^-)=\eta_A(a^+,d^-)+\eta_A(b^+,c^-)$.
\bigbreak

For the order $a^+,b^-,c^+,d^-$, something slightly different occurs.
In the pairings $\{a,b\}$ and $\{c,d\}$, we get the points between $a^+$ and $b^-$ and between $c^+$ and $d^-$ counted exactly once each, and no other
point is counted at all.  In the pairings $\{a,d\}$ and $\{b,c\}$, we get the points  between $a^+$ and $b^-$ and between $c^+$ and $d^-$ counted
exactly twice each, while every other point is  counted once each.  Fortunately, the total over all the letters in $A$ is 0, so the two sums are the
same.
\bigbreak

\section{Invariance under Reidemeister moves}\label{ReidInv}
In this section, we shall take up in turn the three flat Reidemeister moves and show that $p_L$ is unchanged by the move.
The Type I move is easy, since the crossing $e$ that occurs has both $e^+$ and $e^-$ in the same codeword $A$, and they are consecutive on $A$.  It
follows that $\eta(e^+,e^-)$ is 0, because either none or all of the other points on $A$ occur between $e^+$ and $e^-$.  Moreover, for any other pair,
in computing $\eta_A(x^+,y^-)$ either both or neither of $e^+$ and $e^-$ are counted, and so the contribution of $e$ to $\eta_A(x^+,y^-)$ is 0, both
with the crossing $e$ and without.
\bigbreak

Type II moves are equally easy. If the crossings $e$ and $f$ are of $A$ with
itself, then we have $e^\epsilon f^{-\epsilon}$ consecutive in $A$, in this order, and we have either $e^{-\epsilon} f^\epsilon$ or $ f^\epsilon
e^{-\epsilon}$ also occurring consecutively in $A$.  The Type II move removes all the $e$ and $f$ letters from the codeword, so we need to realize that
$\eta_A(e^+,e^-)+\eta_A(f^+,f^-)=0$.  In both cases, the two terms in the sum between them include all the letters of $A$ exactly once, except possibly
none of $e^{+}$, $e^{-}$, $f^+$ and $f^-$, and therefore is 0.  
\bigbreak

In the cases where the two strands are in different components $A$ and $B$, we may assume they are paired,
 since each strand has both a positive and a negative crossing.  In this case, $\eta_A(e^+,f^-)$ and $\eta_B(f^+,e^-)$ are both zero, since they count
either none or all of the other signs of their codeword, and this is obviously independent of the orientations of the strands.
\bigbreak

For the Type III move, we recall that  it suffices to treat the case the strands make a cyclic triangle as in Figure 3
(see, for example, \cite{kauf2}).    We remark that the Type III move inverts within the codeword the order of all three pairs of consecutive letters,
corresponding to the two crossings  in each of the three strands.
\bigbreak

$$ \picill4inby2.5in(Figure3)  $$
\begin{center}
{\sc Figure 3 - A Cyclic Third Flat Reidemeister Move}
\end{center}

Label the crossings  as in Figure 3.  We shall deal with two cases:  in the first, the crossing (we may assume it to be 1)
is of two strands in the same component, while in the second, the crossing is of two strands in different components.  
\bigbreak

For the first case, let $A$ be the component containing the two strands.  
Then $A$'s codeword has $1^+3^-$ and $2^+1^-$.  The Type III move inverts both of these.  Initially, $\eta_A(1^+,1^-)$ counted both the $3^-$ and the
$2^+$, while after the move, it counts neither.  Any other $+$ or $-$ is counted either  by both or by neither.  Thus, $\eta_A(1^+,1^-)$ is unchanged.  
\bigbreak

In the second case, let $A$ be the component containing the $13$-strand, let $B$ be the component containing the $21$-strand,
and let $e$ be the other $AB$-crossing paired with $1$.  We see that from before the Type III move to after it, $\eta_A(1^+,e^-)$ loses the $3^-$ and
$\eta_B(e^+,1^-)$ loses the $2^+$.  Thus, $\eta_A(1^+e^-)+\eta_B(e^+,1^-)$ is unchanged.  This completes the proof that $p_L$ is invariant.
\bigbreak

\section{Relation to filamentations for knots and links}
According to Hrencecin and Kauffman \cite{hren}, a {\em filamentation of a flat virtual knot $K$\/} is a partition
 of the letters of the Gauss codeword into singletons and pairs, so that, for each part $\{x\}$ (the singletons are the {\em monofilaments\/} of the
partition), $\eta(x^+,x^-)=0$, while for each part $\{x,y\}$ (the pairs are the {\em bifilaments\/}) of the partition, $\eta(x^+,y^-)+\eta(y^+,x^-)=0$.
\bigbreak

\noindent {\bf Theorem B.} {\it A flat virtual knot $K$ has a filamentation if and only if $p_K\equiv 0$.}
\bigbreak

\noindent {\bf Proof.} Suppose $\mathcal P$ is a filamentation for $K$.  The coefficient of $t^n$ is just the sum of all
the $\eta(x^+,x^-)$ such that $|\eta(x^+,x^-)|=n$.    For each monofilament $\{x\}$ of $\mathcal P$, we have $\eta(x^+,x^-)=0$, so this contributes 0
to $p_K$.  For each bifilament $\{x,y\}$, we have $\eta(x^+,x^-)=-\eta(y^+,y^-)$, so that, for $n=|\eta(x^+,x^-)|$, these combine to contribute 0 to the
coefficient of $t^n$.  Since $\mathcal P$ is a partition of the chords of $K$, all contributions to $p_L$ are accounted for precisely once, showing
$p_K\equiv 0$.
\bigbreak

The converse is similar.  The monofilaments are the chords $x$ such that $\eta(x^+,x^-)=0$.
For each $n>0$, since $p_K=0$, then number of chords $x$ such that $\eta(x^+,x^-)=n$ is equal to the number of chords $x$ such that
$\eta(x^+,x^-)=-n$.  Thus, these can be partitioned into pairs $\{x,y\}$ so that $\eta(x^+,x^-)=-\eta(y^+,y^-)$.  These pairs are the bifilaments.
$\hfill \Box$
\bigbreak

Schellhorn generalized filamentations to links \cite{schell}.  A {\em filamentation of a flat virtual link $L$\/}  is a partition
of the letters of all the Gauss codewords into singletons and pairs, so that, for each monofilament $\{x\}$, the occurrences of $x^+$ and $x^-$ are in
the same codeword $A$ and $\eta_A(x^+,x^-)=0$, while for each bifilament $\{x,y\}$, $x^+$ and $y^-$ are in the same codeword $A$, while $x^-$ and $y^+$
are in the same codeword $B$ (possibly $A=B$) and $\eta_A(x^+,y^-)+\eta_B(y^+,x^-)=0$.
\bigbreak

In much the same way as the first half of the proof of Theorem B, we obtain the following.
\bigbreak

\noindent {\bf Theorem C.} {\it If a flat virtual link $L$ has  a filamentation, then $p_{L} \equiv 0$.}
\bigbreak

On the other hand, if the link polynomial is 0, there need not be a 
filamentation; we give such an example in the next section.  
Nevertheless, the polynomial being 0 is useful information, because, 
as we are about to show, it will quickly allow us to determine if 
there is a filamentation.
\bigbreak

Let $A$ and $B$ be two components of the link $L$.  In a filamentation 
for $L$, the $AB$-crossings are
partitioned into pairs $\{x,y\}$ so that $A$ has  $x^+$ and $y^-$, 
while $B$ has $x^-$ and $y^+$, and so that 
$\eta_A(x^+,y^-)+\eta_B(y^+,x^-)=0$.  We shall refer to such a 
partition of the $AB$-crossings as a {\em $0$-sum partition\/}.
\bigbreak

\noindent The main technical point is the following.
\bigbreak

\noindent{\bf Lemma D.}  {\em Suppose there is a 0-sum 
partition of the $AB$-crossings.  Suppose $x$ and $y$ are 
$AB$-crossings such that $A$ has  $x^+$ and $y^-$, $B$ has $x^-$ and 
$y^+$, and $\eta_A(x^+,y^-)+\eta_B(y^+,x^-)=0$.  Then there is a 0-sum 
partition of the $AB$-crossings that contains $\{x,y\}$.}
\bigbreak

Given that $p_L\equiv 0$, Lemma D allows us to find a filamentation, if 
it exists, one pair at a time.  To see this, note that, as described in 
the proof of Theorem B, the proposed filamentation will contain the 
filamentations from each component $A$ -- we get the monofilaments 
$\{x\}$ for which $\eta_A(x^+,x^-)=0$ and the bifilaments $\{x,y\}$ 
consisting of a pair of self-crossings of $A$ for which 
$\eta_A(x^+,x^-)=-\eta_A(y^+,y^-)$.
\bigbreak

For two components $A$ and $B$ of the link, we first find a pair 
$\{x,y\}$ of $AB$-crossings such that $A$ has  $x^+$ and $y^-$, $B$ has 
$x^-$ and $y^+$, and $\eta_A(x^+,y^-)+\eta_B(y^+,x^-)=0$.  By Lemma D, 
we can put this pair into the 0-sum partition for the $AB$-crossings, 
remove these two letters from consideration, and repeat.  Either we 
obtain a 0-sum partition for the $AB$-crossings (and repeating this for 
all pairs of components completes the filamentation) or at some stage 
we cannot find another pair, in which case we deduce there is no 
filamentation.
\bigbreak

\noindent{\bf Proof of Lemma D.}    Let $\mathcal P$ be a 0-sum 
partition of the $AB$-crossings.  We may suppose $\{x,y\}\notin 
\mathcal P$, but that $\{x,z\}$ and $\{w,y\}$ are in $\mathcal P$.  
Then $z^-$ and $w^+$ are in $A$.  When we proved that $p_L$ is 
invariant under different partitions, we showed that 
$\eta_A(x^+,z^-)+\eta_A(w^+,y^-)=\eta_A(x^+,y^-)+\eta_A(w^+,z^-)$.  
With $B$ in place of $A$, we also have 
$\eta_B(z^+,x^-)+\eta_B(y^+,w^-)=\eta_B(y^+,x^-)+\eta_B(z^+,w^-)$.  
Summing these and recalling that $\eta_A(x^+,z^-)+\eta_B(z^+,x^-)=0$, 
$\eta_A(w^+,y^-)+\eta_B(y^+,w^-)=0$, and 
$\eta_A(x^+,y^-)+\eta_B(y^+,x^-)=0$, we conclude that 
$\eta_A(w^+,z^-)+\eta_B(z^+,w^-)=0$.  Thus, replacing $\{x,z\}$ and 
$\{w,y\}$ in $\mathcal P$ with $\{x,y\}$ and $\{w,z\}$ yields a new 
0-sum partition of the $AB$-crossings, as required. \hfill $\Box$
\bigbreak

\section{Two examples}\label{examples}

In this section, we provide two simple examples.  The first is  the link $L_1$ of Figure 1,
which has the property that $p_{L_1}=0$,  and yet there is no filamentation.  (Both these claims are very easily checked.)   The latter is what allows
us to conclude that the link is not trivial.  The second, in Figure 4, is an example of a flat virtual link $L_2$ with non-zero polynomial.  As remarked
in the introduction, this implies that these flat universes are non-trivial, and  any virtual link with either of these universes is non-trivial.  
\bigbreak

$$ \picill4inby2in(Figure4)  $$
\begin{center}
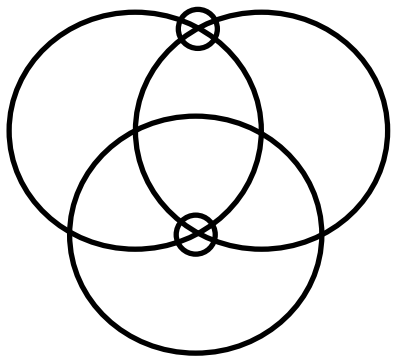
{\sc Figure 4 - A Flat Version of the Borrommean Rings}
\end{center}

\end{document}